\documentclass{article}

\usepackage{arxiv}

\usepackage[utf8]{inputenc} 
\usepackage[T1]{fontenc}    
\usepackage{hyperref}       
\usepackage{url}            
\usepackage{booktabs}       
\usepackage{amsfonts}       
\usepackage{nicefrac}       
\usepackage{microtype}      
\usepackage{lipsum}
\usepackage{graphicx}
\usepackage{amsmath}
\usepackage{algorithm}
\usepackage{algpseudocode}
\usepackage{authblk}

\graphicspath{ {./images/} }

\title{Iterative Belief Propagation for Sparse Combinatorial Optimization}

\author[1]{Sam Reifenstein}
\author[1,2]{Timothée Leleu}
\affil[1]{NTT Research Inc}
\affil[2]{Stanford University}

\begin{document}
\maketitle
\begin{abstract}
In this note we study an iterative belief propagation (IBP) algorithm and demonstrate it's ability to solve sparse combinatorial optimization problems. Similar to simulated annealing (SA) \cite{Kirkpatrick1983_SA}, our IBP algorithm attempts to sample from the Boltzmann distribution of the objective function but also uses belief propagation (BP) \cite{Braunstein2005_SP, Pearl1982_BP} to improve convergence.
\end{abstract}


\section{Introduction and Algorithm Overview}
Combinatorial optimization problems are ubiquitous in applications across many field of science and engineering. However, because many of these problems are NP-Hard, they are expected to take exponential time to solve in worst case by exact algorithms. Thus, when exact algorithms do not work, heuristic approaches such as simulated annealing (SA) provides a more practical method for finding exact and approximate solutions \cite{Kirkpatrick1983_SA}. 
\\
\\
SA works by attempting to for sampling from the Boltzmann distribution of the optimization problem. That is, given a solution a cost function $C(x)$, where $x \in X$ is a possible solution, and ``inverse temperature" $\beta$, we want to sample from the following probability distribution:
\begin{equation}
    P(x) = \frac{1}{Z}e^{-\beta C(x)}.
\end{equation}
This is useful for optimization because in the limit that $\beta$ becomes large, $P(\text{argmin}_{X}) \approx 1$ (assuming there is one minimizer). SA generally uses the Metropolis-Hastings criterion to iteratively update a solution $x$, aiming to sample from this distribution while gradually increasing (``annealing") the inverse temperature $\beta$ from a low to a high value. This update step can be interpreted as choosing a random sub-problem of the full problem (of just one variable in this case) and sampling from the Boltzmann distribution of that sub-problem while keeping all other variables fixed. SA has the useful property that it will, given sufficient iterations, converge to the correct solution with nonzero probability for any problem. However, in the worst case, this convergence may require exponential time with respect to the problem size.
\\
\\

Belief propagation (BP) is an algorithm which allows to approximately sample from a high dimensional discrete probability distribution when the probability function can be expressed as a product of factors that involve a small subset of the variables. If the graph that connects factors and variables forms a tree, the BP algorithm is know to converge exactly in polynomial time\cite{Pearl1982_BP}. In the context of optimization we can imagine that we have an objective function which can be expressed as the sum of smaller objectives that involve small subsets of variables. Once we exponentiate this objective to get the Boltzmann distribution these small objectives turn into factor and if the connectivity forms a tree then this optimization problem will be solvable in polynomial time using BP. To make things more straightforward, from now on we will primarily discuss the NP-hard QUBO (quadratic unconstrained binary optimization) problem even though the algorithms that we outline are more general. The QUBO problem is defined as follows. We want to minimize a set of binary variables $x \in \{0,1\}^{N}$ over a quadratic cost function defined as
\begin{equation}
    C(x) = \sum_{i,j} Q_{ij} x_i x_j
\end{equation}
where $Q$ is a symmetric $N\times N$ ``QUBO matrix". In the context of BP, each nonzero off-diagonal entry of $Q$ will correspond to a factor connected to two variables while the diagonal elements of $Q$ will become factors connected to just one variable. So, if $Q$ corresponds to a graph that is acyclic then the QUBO can be solved in polynomial time with BP. Belief propagation is also known to work in some cases when the graph has loops (sometimes called ``loopy belief propagation"), but in general it does not converge reliably to the correct marginals \cite{Braunstein2005_SP}. Thus, unlike SA, BP is not typically useful for solving hard optimization problems.
\\
\\
We study an iterative belief propagation (IBP) algorithm that attempts to bridge the gap between BP and SA as follows (not to be confused the iterative belief propagation algorithm of \cite{Dechter2012_IBP}). Given a QUBO graph with sparse connectivity, it is possible to choose a subset of variables such that the subgraph formed by them is a tree. If the given QUBO graph is sparse and more tree-like, these sub-graphs will be larger, whereas if the graph has full connectivity these subgraphs will have size 2 at most. IBP works by randomly choosing one of these subgraphs (which can be done efficiently) and applying BP on it. This allows us to sample from the Boltzmann distribution for a subset of variables similar to SA, except unlike SA many variables will be updated at once. Iteratively choosing these sub-trees and updating $x$ using BP will allow us to construct a Markov chain which converges to the Boltzmann distribution like SA. The nice thing about IBP is that in the limiting case of a QUBO that is already a tree, the algorithm will identical to regular BP and thus converge in polynomial time (given $\beta$ is chosen correctly). On the other hand, if the problem is fully connected IBP will reduce to a slightly modified version of SA in which two variables are updated at once instead of one. Many problems of industrial and academic interest lie somewhere in-between these two extremes thus it is reasonable to postulate that IBP can outperform both BP and SA in some practically relevant cases.

\section{Numerical Results}
\begin{figure}[hb]
    \centering
    \includegraphics[width=0.3\linewidth]{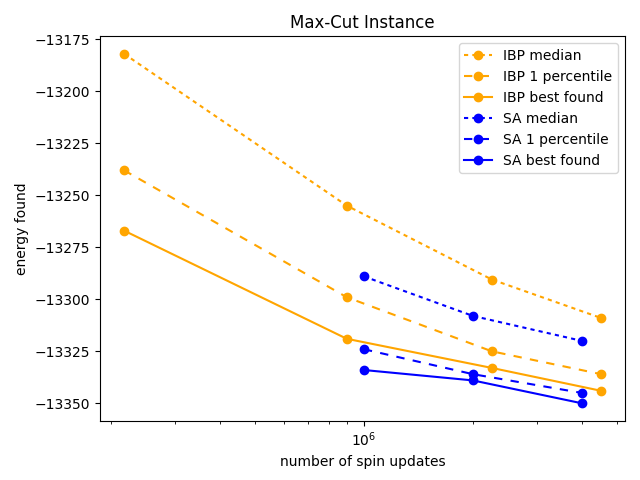}
    \includegraphics[width=0.3\linewidth]{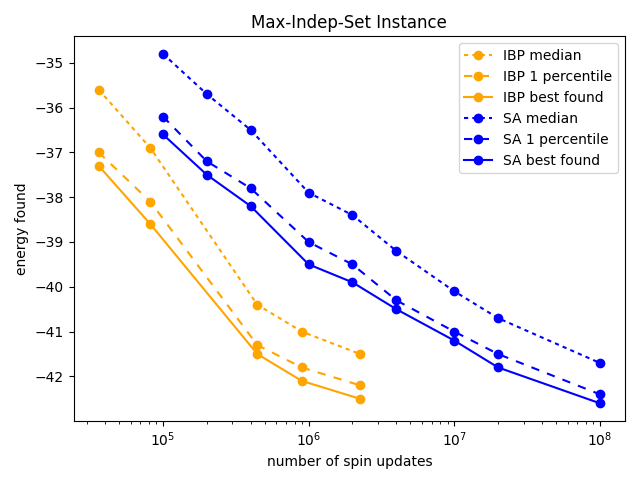}
    \includegraphics[width=0.3\linewidth]{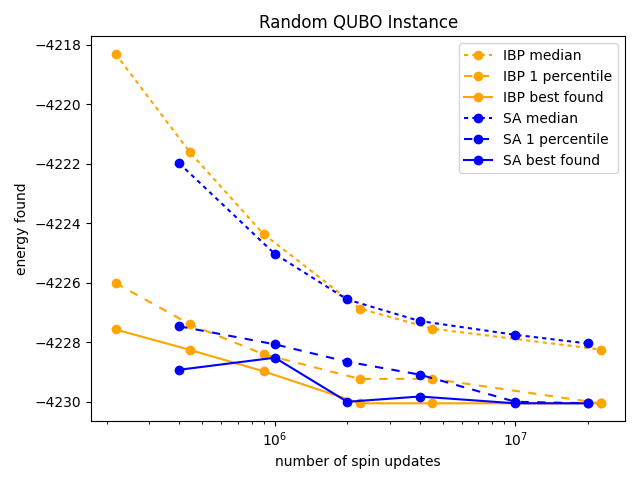}
    \caption{Median, 1st percentile and best objective values found by SA and IBP with respect to total spin updates on Left: a Max-Cut instance, Middle: a Maximum Independent Set instance and Right: a sparse random QUBO instance. All three instance are generated from Erdős–Rényi random grpahs of size $N=2000$ and density $1\%$. The number of spin updates for IBP is defined as the total size of all sub-trees considered (average size pf sube-tree multiplied by number of iterations) while the number of spin update for SA is number of sweeps muptiplied by $N$. The same geometric annealing sechdule is used for $\beta$ for both algorithms.}
    \label{fig:random_problem_results}
\end{figure}

In this section we will briefly outline some numerical results in which IBP is compared to SA on a few sparse QUBO instances. We consider the three following random classes of QUBO instances:
\begin{itemize}
    \item \textbf{Max-Cut:} Max-Cut problem of given random Erdős–Rényi graph is converted to QUBO form (all edges have weight 1)
    \item \textbf{Maximum Independant Set:}  Maximum Independant Set problem of random Erdős–Rényi graph is converted to QUBO form 
    \item \textbf{Random Sparse QUBO:} QUBO Matrix elements corresponding to edges of graph are set to a uniform random number in the range $[-1,1]$ while the rest are set to 0. Diagonal QUBO elements are also set to a uniform random number in $[-1,1]$.
\end{itemize}
In figure \ref{fig:random_problem_results} we show results for IBP and SA on one random instance of size $N=2000$ and density $1\%$ from each of the three classes. We find that difference between IBP and SA varies greatly across different problem types, with IBP performing the best relative to SA on the Maximum Independent Set problem and worst on the Max-Cut problem. In the best case IBP reduces the numbers of spin updates required by orders of magnitude.

\section{Algorithm Details}
\subsection{Algorithm Pseudo-code}
The IBP algorithm takes in QUBO matrix $Q$, number of replicas $R$, and list of inverse temperatures $\text{BETAS}$ and outputs $x$ a list of $R$ possible solutions to the QUBO problem. $D$ denotes the maximum degree of the graph and $D_{sup}$ denotes the maximum degree of the subtree.
\begin{algorithm}[H]
    \caption{IBP \textbf{Input:} QUBO matrix $Q$, number of replicas $R$, and list of inverse temperatures $\text{BETAS}$ }
    \begin{algorithmic}[1]
    
    \State $R$ copies of $x$ initialized randomly
    \For{$\beta$ in $\text{BETAS}$}
        \State \textit{Select} random sub-tree of size $M$ (same for all replicas) \Comment{see section \ref{sec: sub-tree}}
        \State \textit{Compute} the off diagonal couplings for the sub-problem (same for all replicas) and diagonal couplings (different for each replica) \Comment{Computational complexity $O(MD)$}
        \State \textit{Perform} BP on sub-problems (in parallel) \Comment{see sections \ref{sec: bp-1},\ref{sec: bp-2}, computational complexity $O(MD_{\text{sub}}^2)$}
        \State \textit{Sample} using BP Soln and update subset of $x$ (in parallel) \Comment{see section \ref{sec: bp-sample}, computational complexity $O(MD_{\text{sub}}^2)$}
    \EndFor
    \State \textit{Return} $x$ as a set of $R$ solution vectors, or best solution found
    
    \end{algorithmic}
\end{algorithm}
\subsection{Selecting Sub-trees}\label{sec: sub-tree}
To select a sub-tree we start by choosing a random variable. Then we can iteratively add variables to the tree that are connected to exactly one member of the current set. This makes it impossible for loops to form on this subset. Once it is impossible to select variables with exactly one connection we stop. There are many different ways to select sub-trees but this simple algorithm seems to work fine for the purposes of IBP. Additionally it is worth noting that when the number of replicas $R$ is large the sub-tree selection will not contribute significantly to the computation time since the same sub-tree is used for all replicas.

\subsection{Belief Propagation}\label{sec: bp-1}
For QUBO, the belief propagation algorithm can be described as follows. For each ordered pair of variables $i,j$ such that $Q_{ij} \neq 0$ BP stores two ``messages" $\mu^{a}_{ij}$ with $a \in \{0,1\}$. These messages are real numbers in the range $[0,1]$ and are typically initialized randomly. The messages are updated iteratively according to:
\begin{equation}\label{eq: bp1}
    (\mu^{a}_{ij})^{*} = \prod_{k \in \mathcal{N}(i), k \neq j} \left[\sum_{b \in \{0,1\}}e^{-\beta Q_{ik} ab}\prod_{l \in \mathcal{N}(k), l \neq j} \mu^{b}_{lk} \right]
\end{equation}
\begin{equation}\label{eq: bp2}
    \mu^{a}_{ij} \rightarrow \frac{(\mu^{a}_{ij})^{*}}{\sum_{b \in \{0,1\}} (\mu^{b}_{ij})^{*}} 
\end{equation}
The set $\mathcal{N}(i)$ refers to all of the connected QUBO variables (which can include $i$ itself). The second equation is not necessary but helps with numerical stability of the algorithm. Once the algorithm has converged the marginal probabilities can be computed as:
\begin{equation}
    (p^{a}_{i})^{*} = \prod_{j \in \mathcal{N}(i)}{\mu^{a}_{ij}} \quad \quad P(x_i = a) = \frac{(p^{a}_{i})^{*}}{\sum_{b\in\{0,1\}}(p^{b}_{i})^{*}}
\end{equation}

\subsection{Logarithmic form of BP}\label{sec: bp-2}
When $\beta$ is large, \eqref{eq: bp1} and \eqref{eq: bp2} involve many products of exponentials that can range over orders of magnitude in values. This poses a problem for typical 64-bit floating point arithmetic. Luckily, the equations can be re-parameterized using logarithms to get rid of this problem. We introduce the new variable $z_{ij} \in \mathbb{R}$ such that:
\begin{equation}
    \mu^a_{ij} = \frac{1}{1 + e^{(1-2a)z_{ij}}}
\end{equation}
(or equivalently $z_{ij} = \log(\mu^1_{ij}/\mu^0_{ij})$). Then we can write the update equations as:
\begin{equation}
    z_{ij} \rightarrow  \sum_{a \in \{0,1\}} (1 - 2a) \sum_{k \in \mathcal{N}(i), k \neq j} \log \left[\sum_{b \in \{0,1\}}\exp \left[{-\beta Q_{ik} ab} + \sum_{l \in \mathcal{N}(k), l \neq j}\log \left(\frac{1}{1 + e^{(1-2b)z_{lk}}} \right) \right] \right]
\end{equation}
\begin{equation} 
 = \sum_{a \in \{0,1\}} (1 - 2a) \sum_{k \in \mathcal{N}(i), k \neq j} \log \left[1 + \exp \left[{-\beta Q_{ik} a} + \sum_{l \in \mathcal{N}(k), l \neq j} z_{lk} \right] \right]
\end{equation}
\begin{equation} \label{eq: bp-final}
 =  \sum_{k \in \mathcal{N}(i), k \neq j} \log \left[ \frac{1 + \exp \left[\sum_{l \in \mathcal{N}(k), l \neq j} z_{lk} \right] }{  1 + \exp \left[{-\beta Q_{ik} a} + \sum_{l \in \mathcal{N}(k), l \neq j} z_{lk} \right]} \right]
\end{equation}
This expression can be further simplified to a form in which numerical overflow errors are impossible.
\subsection{Sampling from the Boltzmann Distribution with BP}\label{sec: bp-sample}
Although BP is typically described as an algorithm to compute marginal probabilities, to actually sample from the desired probability distribution requires a few more steps. To start off, we can pick a random variable and choose its value based on the computed marginal probability. However, to choose the next variable we need to use a conditional probability on our original choice. Importantly, for connected variables these conditional probabilities can also be computed using the messages as:
\begin{equation} \label{eq: cond-samp}
    (P^{ab}_{ij})^{*} = e^{-\beta ab Q_{ij}} \prod_{k \in N(i)} \prod_{l \in \mathcal{N}(k), l \neq j} \mu^{a}_{lk} \quad \quad P(x_i = a \mid x_j  = b) = \frac{(P^{ab}_{ij})^{*}}{ \sum_{c \in \{0,1\}} (P^{cb}_{ij})^{*} }
\end{equation}
Using this formula we can propagate along the tree starting with the initial sampled variable and sample all variables of the sub-tree.
\subsection{Computational Complexity}
When the number of replicas is large, the only contributions to the computation time are the operations that need to be done in parralel for all replicas. This comes down to 3 operations, calculation of the diagonal (self-interacting) elements for the sub-problem (line 5 of the pseudo-code), computing the BP iterations, and sampling the new state. The first operation takes roughly $O(MD)$ time where $D$ is the degree of the graph since roughly $MD$ QUBO elements need to be summed. For BP and sampling procedures require at least $O(MD_{sub}^2)$ because of the nested sums in equations \eqref{eq: bp-final, eq: cond-samp}. However, with a better implementation this could likely be reduced to $O(M)$ sums (because a tree has same number of nodes as edges). So, as described in this document the algorithm takes roughly $O(MD + C MD_{\text{sub}}^2)$. Per spin update this becomes $O(D + C D_{sub}^2)$ which is close to $O(D)$ complexity of SA. With an efficient implementation it is likely that the complexity per spin update can be very close to that of SA in most cases. However the implementation used for the numerical results presented here is much slower per spin update than currently available SA implementations in C.

\section{Conclusion and Outlook}
In conclusion we propose iterative belief propagation, a generic algorithm for sparse combinatorial optimization. We have shown numerically that in certain cases IBP is able to outperform its close relative, simulated annealing. Future works will consider a much more complete benchmark of this algorithm and also attempt to clarify exactly when IBP can provide improvement over SA and other state of the art algorithms. Additionally, although this study was limited to QUBO, future works will consider the performance of this algorithm on other types of NP-Hard combinatorial optimization problems such as boolean satisfiability.
\bibliographystyle{unsrt}  
\bibliography{references}

\end{document}